\documentclass{amsart}

\usepackage{amsmath}
\usepackage{amssymb}
\usepackage{latexsym}
\usepackage{enumerate}
\usepackage{amscd}
\usepackage[all]{xy}
\usepackage{graphics}

\newcommand{\HH}{\mathrm{H}}                       

\newcommand{\Z}{\mathbf{Z}}

\newcommand{\mac}{\mathcal{M}_{g,1}}                

\newtheorem{lemma}{Lemma} 
\newtheorem{proposition}{Proposition} 
\newtheorem{theorem}{Theorem} 
\newtheorem{corollary}{Corollary}

\begin{document}

\title{Trivial Cocycles and Invariants of Homology $3$-Spheres}

\author{Wolfgang Pitsch}
\address{Wolfgang Pitsch \\ Universitat Aut\`onoma de Barcelona \\ Departament de Ma\-te\-m\`ati\-ques\\
E-08193 Bellaterra, Spain}
\email{pitsch@mat.uab.es}

\thanks{The author is supported  by MEC grant MTM2004-06686 and  by the program Ram\'on y Cajal, MEC, Spain}

\begin{abstract}
We study the relationship between trivial cocycles on the Torelli group and invariants of oriented integral homology $3$-spheres. We apply this study to give a new purely algebraic construction of the Casson invariant and prove in this setting its surgery properties. As a by-product we get a new $2$-torsion cohomology class in the second integral cohomology of the Torelli group.
\end{abstract}


\maketitle

\section*{Introduction}\label{sec Intro}
If one tries to understand $3$-manifolds by  ``cut and paste" techniques one faces two different paths: either one can concentrate the difficulties in the pieces and have ``simple'' glueing maps (see Kneser's prime de\-com\-po\-si\-tion or Thurston de\-com\-po\-si\-tion into geometric pieces) or one can concentrate the difficulty into the glueing maps and get ``simple" pieces. In the latter path one finds the theory of Heegaard splittings, the pieces are handlebodies and the glueing problems are encompassed within the mapping class groups of oriented surfaces, $\mathcal{M}_{g,1}$. It is natural then to try to construct invariants of $3$-manifolds out of the algebraic properties of this groups. For general $3$-manifolds this strategy has been adopted for instance by Birman in \cite{MR0375318} but is usually hopelessly difficult, for the structure of $\mathcal{M}_{g,1}$ is quite involved. In this paper we will concentrate on the subclass of integral homology $3$-spheres. In this case the group that controlls the glueings is the Torelli group, $\mathcal{T}_{g,1}$. This approach has been that of Morita in \cite{Mo1} where he constructed a function out of the Torelli group that he showed to coincide point-wise with the Casson invariant. Although this did not succeed into a new construction of this important invariant it was the starting point for a fruitfull exploration the interplay between the Casson invariant and algebraic properties of the Torelli group \cite{Mo1, Mo2, MR1133875}.

In this paper we will give a general framework to construct invariants of ho\-mo\-logy spheres in a purely algebraic setting. We will show that the algebraic problems  boil down to a low-dimensional cohomological problem. As an example we will give a construction of an invariant of homology spheres and by proving the ``surgery formulas" we will show that it coincides with the Casson invariant.

Denote by $\mathcal{V}(3)$ the set of diffeomorphism classes of compact, closed and oriented smooth $3$-manifolds and by $\mathcal{S}(3) \subset \mathcal{V}(3)$ the subset of homology spheres, that is diffeomorphisms classes that have the same integral homology as the standard $3$-sphere $\mathbf{S}^3$.
Let $\Sigma_g$ denote an oriented surface of genus $g$ standardly embedded in the oriented $3$-sphere $\mathbf{S}^3$. In particular $\Sigma_g$ separates $\mathbf{S}^3$ into two genus $g$ handlebodies $\mathbf{S}^3 = \mathcal{H}_g \cup - \mathcal{H}_g$ with opposite induced orientation. Denote by $\mathcal{M}_{g,1}$ the mapping class group of $\Sigma_g$, that is the group of  orientation-preserving diffeomorphisms of $\Sigma_g$ which are the identity on a small fixed disc modulo isotopies which fix that small disc pointwise. The embedding $\Sigma_g \hookrightarrow \mathbf{S}^3$ determines three natural subgroups of $\mac$, namely the subgroup $\mathcal{B}_{g,1}$ of mapping classes that are restrictions of diffeomorphisms of the first handlebody $\mathcal{H}_g$, the subgroup $\mathcal{A}_{g,1}$ of mapping classes that are restrictions of diffeomorphisms of the second handlebody $-\mathcal{H}_g$ and their intersection $\mathcal{AB}_{g,1}$.

From the theory of Heegaard splittings we learn that any element in $\mathcal{V}(3)$ can be obtained by cutting $\mathbf{S}^3$ along $\Sigma_g$ for some $g$ and glueing back the two handlebodies by some element $\phi \in \mathcal{M}_{g,1}$. The lack of injectivity of this construction is controlled by the subgroups $\mathcal{B}_{g,1}$ and $\mathcal{A}_{g,1}$. More precisely there is a natural injective sta\-bili\-zation map $\mathcal{M}_{g,1} \hookrightarrow \mathcal{M}_{g+1,1}$, which is compatible with the definitions of the above subgroups and one gets a well-defined bijective map~:
\[
\begin{array}{ccc}
\displaystyle\lim_{g \to \infty} \mathcal{B}_{g,1} \backslash \mac \slash \mathcal{A}_{g,1} & \stackrel{\sim}{\longrightarrow} & \mathcal{V}(3) \\
\phi & \longmapsto & \mathbf{S}^3_\phi = \mathcal{H}_g \cup_\phi -\mathcal{H}_g.
\end{array}
\]

Thus any problem on 3-dimensional manifolds can be translated into a problem on the mapping class group. In particular any invariant $F~: \mathcal{V}(3) \rightarrow \Z$ can be viewed as a compatible family of functions on the mapping class groups $\mathcal{M}_{g,1}$ which are constant on double cosets.

If we restrict our study to $\mathcal{S}(3)$ the situation becomes more tractable. First we can restrict our attention to those mapping classes that act trivially on the homology of the underlying surface. Recall that the Torelli group $\mathcal{T}_{g,1}$ is defined as  the kernel of the natural map $\mathcal{M}_{g,1} \longrightarrow \text{Aut } (H_1(\Sigma_g;Z))$. The above bijection induces a new bijection \cite{Mo1}~:
\[
\begin{array}{ccc}
\displaystyle\lim_{g \to \infty} \mathcal{B}_{g,1} \backslash \mathcal{T}_{g,1} \slash \mathcal{A}_{g,1} & \stackrel{\sim}{\longrightarrow} & \mathcal{S}(3) \\
\phi & \longmapsto & \mathbf{S}^3_\phi = \mathcal{H}_g \cup_{\phi} -\mathcal{H}_g,
\end{array}
\]
where $\mathcal{B}_{g,1} \backslash \mathcal{T}_{g,1} \slash \mathcal{A}_{g,1}$ stands for those mapping classes in $\mathcal{M}_{g,1}$ that contain an element of the Torelli group. Denote by $\mathcal{TB}_{g,1}$ (resp. $\mathcal{TA}_{g,1}$) the group $\mathcal{T}_{g,1} \cap \mathcal{B}_{g,1}$ (resp. $\mathcal{T}_{g,1} \cap \mathcal{A}_{g,1}$). The induced equivalence relation on the Torelli group has an intrinsic description:

\begin{theorem}\label{thm intro1}
Two elements $\phi, \psi \in\mathcal{T}_{g,1}$ are in the same double coset in $\mathcal{B}_{g,1} \backslash \mac \slash \mathcal{A}_{g,1}$ if and only if there exist maps $\xi_b \in \mathcal{TB}_{g,1}, \xi_a \in \mathcal{TA}_{g,1}$ and $\mu \in \mathcal{AB}_{g,1}$ such that
\[
\phi = \mu \xi_b  \psi \xi_a \mu^{-1}.
\]1
\end{theorem}

The conjugacy part of this equivalence relation is the key tool of our study. Consider an integral valuated invariant of homology spheres $F~: \mathcal{S}(3) \rightarrow \Z$. By the above bijection and Theorem we can view $F$ as a family of compatible  functions $F_g$ (i.e $F_{g+1} \vert_{\mathcal{T}_g,1} = F_g)$ that are constant on the double coset  classes $\mathcal{TB}_{g,1} \backslash \mathcal{T}_{g,1} \slash \mathcal{TA}_{g,1}$ and invariant under conjugation by $\mathcal {AB}_{g,1}$. To any such family of functions we associate a family of \emph{trivialized} $2$-cocycles on the Torelli groups $C_g(\phi,\psi) = F_g(\phi) + F_g(\psi) - F_g(\phi \psi)$. It turns out that these functions are not trivial unless $F_g$ is itself trivial. Since $\mathcal{T}_{g,1}$ is not perfect there is a difference between trivialized cocycles and trivial cocycles. One might wonder what conditions we should impose on a family $(C_g)$ of trivial $2$-cocycles on the Torelli groups such that from their trivializations one can extract a family of compatible trivializations $(F_g)$ that reassemble into an invariant of homology spheres $F~: \mathcal{S}(3) \rightarrow \Z$. Notice that the maps $F_g$ are necessarily $\mathcal{AB}_{g,1}$-invariant trivializations of the cocycles.

The cocycles $C_g$ inherit the following properties of the maps $F_g$:
\begin{enumerate}
\item[(1)] The cocycles $C_g$ are be compatible $C_{g+1} \vert_{\mathcal{T}_{g,1} \times \mathcal{T}_{g,1} } = C_g$.
\item[(2)] The cocycles $C_g$ are $0$ on $\mathcal{TB}_{g,1} \times \mathcal{T}_{g,1} \cup \mathcal{T}_{g,1} \times \mathcal{TA}_{g,1}$.
\item[(3)] The cocycles $C_g$ are invariant under conjugation by $\mathcal{AB}_{g,1}$.
\end{enumerate}

Then, the existence of an $\mathcal{AB}_{g,1}$-invariant trivialization of the cocycle $C_g$ is controlled by a cohomology class, the torsor~:
\[
\rho(C_g) \in \mathrm{H}^1(\mathcal{AB}_{g,1}; \wedge^3 \mathrm{H}_1(\Sigma_g; \mathbf{Z})).
\]

The three conditions above and the nullity of the torsor turn out to be not only necessary but also sufficient~:

\begin{theorem}\label{thm intro2}
A family of cocycles $(C_g)_{g \geq 3}$ on the Torelli groups $\mathcal{T}_{g,1}$, $ g \geq 3$, satisfying conditions $(1)-(3)$ provides a compatible family of trivializations $F_g : \mathcal{T}_{g,1} \rightarrow \Z$ that reassemble into an invariant of homology spheres
\[
\displaystyle \lim_{g \to \infty} F_g : \mathcal{S}(3) \rightarrow \Z
\]
if and only if the following two conditions hold:
\begin{enumerate}
\item[(i)] The associated cohomology classes $[C_g] \in \mathrm{H}^2(\mathcal{T}_{g,1}; \Z)$ are trivial.
\item[(ii)] The associated torsors $\rho(C_g) \in \mathrm{H}^1(\mathcal{AB}_{g,1}; \wedge^3 \mathrm{H}_1(\Sigma_g))$ are trivial.
\end{enumerate}
In this case the maps $F_g$ are the unique $\mathcal{AB}_{g,1}$-invariant trivializations of the cocycles $C_g$.
\end{theorem}

Obviously, constructing (trivial) $2$-cocycles directly on the Torelli group is still a difficult problem but instead one could try to pull-back known $2$-cocycles defined on homomorphic images of the Torelli group. We successfully apply this strategy to the the Johnson homomorphism $\tau~: \mathcal{T}_{g,1} \rightarrow \wedge^3 \mathrm{H}_1(\Sigma_g; \Z).$ 

\begin{theorem}\label{thm intro3}
The unique $2$-cocycles on $\mathrm{H}_1(\Sigma_g; \Z)$ whose pull-back along the Johnson ho\-mo\-mor\-phism satisfy conditions $(1)-(3)$ are of the form $nJ_g$, $n \in \Z$ for an explicit $2$-cocycle $J_g$. Moreover:
\begin{enumerate}
\item[(1)] The pull-backs of the cocycles $2J_g$ and the associated torsors $\rho(2J_g)$  are trivial.

\item[(2)] The associated invariant is equal to the Casson invariant. 
\end{enumerate}
\end{theorem}

Moreover~:

\begin{theorem}\label{thm intro4}
\begin{enumerate}
\item[1)] The pull-backs of the cocycles $J_g$ on the Torelli groups are not trivial and define stable $2$-torsion cohomology classes $[J_g] \in \mathrm{H}^2( \mathcal{T}_{g,1}; \mathbf{Z}).$ 
\item[2)] Vieweing  the Rohlin invariant as a familly of  classes $R_g \in \mathrm{H}^1(\mathcal{T}_{g,1}; \Z/2\Z)$, we have
\[
\beta_\Z  (R_g) = [J_g],
\]
where $\beta_\Z$ stands for the integral Bockstein operation.
\end{enumerate}
\end{theorem}

We wish to point out that this construction of an invariant as a compatible familly of trivializations of the pull-backs of cocycles $2J_g$ is independent of the construction of the Casson invariant.

Indeed we prove directly from the algebraic relations
\[
F_g(\phi) + F_g(\psi) - F_g(\phi \psi) = 2J_g(\tau(\phi), \tau(\psi))
\]
satisfied by the maps $F_g$ that our invariant $F$ satisfies also the following surgery formulas  properties (see Propositions \ref{prop surgery} and \ref{prop boundarylink} in Section \ref{sec Casson})~:
\begin{proposition}\label{prop intro}
Let $K\subset M$ be a knot in the oriented homology sphere $M$. For an integer $n \geq 1$ denote by $K_{n}$ the result of performing a $\frac{1}{n}$-Dehn surgery on $K$. 
\begin{enumerate}
\item The difference $F'(K) = F(K_{n+1}) -F(K_{n})$ is independent of $n$.
\item If the knots $K$ and $L$ bound disjoint Seifert surfaces in $M$ then the alter\-nating sum
\[
F''(K,L) = F(K_{k+1},L_{l+1}) - F(K_k,L_{l+1}) - F(K_{k+1},L_l) + F(K_k,L_l)
\]
is independent of the integers  $k$ and $l$ and in fact is $0$.
\end{enumerate}
\end{proposition}

As was proved by Casson an invariant of homology spheres that satisfies the above two properties and which vanishes on $\mathbf{S}^3$ is proportional to ``the'' Casson invariant.

Here is the plan of this paper. In Section \ref{sec Heegardsplit} we turn back to the definition of the groups $\mathcal{A}_{g,1}$, $\mathcal{B}_{g,1}$, $\mathcal{AB}_{g,1}$, we describe their actions on the first homology and homotopy groups of the underlying surface and we prove Theorem \ref{thm intro1}. In Section \ref{sec Trivcocycle} we study the relationship between trivial cocycles on the Torelli groups and invariants of homology spheres. In particular we prove Theorem \ref{thm intro2} up to a technical Lemma which is delayed until Section \ref{sec generatosLT}. In Section \ref{sec Casson} we apply our results to give a purely algebraic construction of the Casson invariant and we prove Theorems \ref{thm intro3} and \ref{thm intro4} and Proposition \ref{prop intro}. Finally, in Section \ref{sec generatosLT} we cope with the proof of the technical Lemma.

\subsection*{General conventions}

The properties of the genus $1$ and $2$ mapping class groups and their subgroups is very peculiar. Since the injectivity of the stabilization map $\mac \hookrightarrow \mathcal{M}_{g+1,1}$ implies that in our case  it is enough to consider large enough values of $g$ we will allways assume that $g \geq 3$. All invariants considered will take the value $0$ on the standard oriented sphere $\mathbf{S}^3$. If $\gamma$ denotes a simple closed curve on the surface $\Sigma_g$, we will denote by $T_\gamma$ the right-hand Dehn twist about $\gamma$.

\section{Heegaard splittings of homology spheres}\label{sec Heegardsplit}

\subsection{The mapping class group and some of its subgroups}\label{subsec mapandsubgr}

\begin{figure}[hbt]
\input{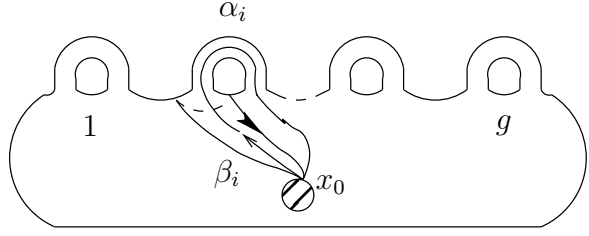}
\caption{Model of $\Sigma_g$}\label{fig model}
\end{figure}

For convenience we fix a model of our genus $g$ surface $\Sigma_g$ as in Figure \ref{fig model}. We denote by $\Sigma_{g,1}$ the complement of the interior of a small disc embedded in $\Sigma_g$. We fix a base point on the boundary of $\Sigma_{g,1}$. The (isotopy classe of) the curves $\alpha_i, \beta_i$ $1 \leq i \leq g$ are free generators of the free group
$ \pi_1 (\Sigma_{g,1},x_0)$. The first homology group of the surface $\mathrm{H}_1 (\Sigma_g; \mathbf{Z}) = H$ is endowed via Poincar\'e duality with a natural symplectic intersection form $\omega: \wedge^2 H \rightarrow \Z$. The homology classes $a_i,b_i$ of the above curves freely generate the abelian group $H \simeq \Z^{2g}$ and define two transverse Lagrangians $A$ and $B$ in $H$.

Denote by $\text{Diff}^+ (\Sigma_g , \text{rel.} D^2)$ the group of orientation preserving dif\-feo\-mor\-phisms of $\Sigma_g$ that are the identity on the fixed small disc, endowed with the compact-open topology. The mapping class group $\mac$ is  the group of connected components $\mac = \pi_0(\text{Diff}^+ (\Sigma_g , \text{rel.} D^2)$.

The natural action of the mapping class group $\mac$ on $H$ clearly preserves the intersection form and we have a short exact sequence where the kernel is known as the Torelli group $\mathcal{T}_{g,1}$:
\[
\xymatrix{ 1 \ar[r] & \mathcal{T}_{g,1} \ar[r] & \mac \ar[r] & \text{Sp}\omega \ar[r] & 1.}
\]

Recall that our surface is standardly embedded in the oriented $3$-sphere $\mathbf{S}^3$. As such it determines two embedded handlebodies $\mathbf{S}^3 = \mathcal{H}_g \cup -\mathcal{H}_g$. By the \emph{inner handlebody} $\mathcal{H}_g$ we will mean the one that is visible in Figure \ref{fig model} and by the \emph{outer handlebody} $-\mathcal{H}_g$ we will mean the complementary handlebody. They are naturally pointed by $x_0 \in \mathcal{H}_g \cap -\mathcal{H}_g$.
 
From these we get three natural subgroups of $\mac$. First, the subgroup of  those mapping classes that are restrictions of diffeomorphisms of the inner handlebody $\mathcal{H}_g$ which we call $\mathcal{B}_{g,1} \subset \mac$. Second, the subgroup of those that are restrictions of the outer handlebody $\mathcal{A}_{g,1} \subset \mac$. Finally, their intersection $\mathcal{AB}_{g,1}$, which may be identified to subgroup of mapping classes that are restrictions of diffeomorphisms of the whole sphere that leave our embedded surface invariant. We denote the groups $\mathcal{T}_{g,1} \cap \mathcal{A}_{g,1}, \mathcal{T}_{g,1} \cap \mathcal{B}_{g,1}$ and $ \mathcal{T}_{g,1} \cap \mathcal{AB}_{g,1}$ respectively by $\mathcal{TA}_{g,1}, \mathcal{TB}_{g,1}$ and $\mathcal{TAB}_{g,1}.$

{\bf Remark:} In this article we deal mostly with mapping class groups relative to a boundary component. Most references, in particular those dealing with the subgroups $\mathcal{A}_{g,1}, \mathcal{B}_{g,1},\mathcal{AB}_{g,1}$ \cite{Gri}, \cite{MR0500977}, \cite{Suz} are written for closed surfaces but lifting their results to the boundary case is not difficult. Indeed, taking isotopy class of diffeomorphisms of the closed surface relative to the base point $x_0$ gives us another mapping class group usually denoted by $\mathcal{M}_{g,\ast}$. The natural ``forgetfull" operation induces a surjective map $\mac \rightarrow \mathcal{M}_{g,\ast}$. Its kernel is ge\-ne\-ra\-ted by a Dehn twist around a curve parallell to the boundary and we get a short exact sequence~:
\[
\xymatrix{1 \ar[r] & \Z \ar[r] & \mac \ar[r] & \mathcal{M}_{g,\ast} \ar[r] & 1. }
\]

The mapping class group of the ``inner" handlebody $\mathcal{H}_g$ relative to the base point can be identified with a subgroup $\mathcal{B}_{g,\ast} \subset \mathcal{M}_{g,\ast}$, where the inclusion is induced by restricting mapping classes to the boundary. Since the aforementioned Dehn twist extends naturally to the handlebody $\mathcal{H}_g$ the preimage of $\mathcal{B}_{g,1} \subset \mac$ of $\mathcal{B}_{g,\ast}$ can be identified as the mapping class group of the handlebody $\mathcal{H}_g$ relative to a small ball $B^3$ such that $B^3 \cap \Sigma_g$ is our distinguished small disk $D^2$ and we get a commutative diagram with vertical arrows induced by restricting to the boundary.
\[
\xymatrix{1 \ar[r] & \Z \ar[r] & \mac \ar[r] & \mathcal{M}_{g,\ast} \ar[r] & 1 \\
1 \ar[r] & \Z \ar[r] \ar@{=}[u] & \mathcal{B}_{g,1} \ar[r] \ar@{^{(}->}[u] & \mathcal{B}_{g,\ast} \ar[r] \ar@{^{(}->}[u] & 1.
}
\]
Similar identifications hold for the groups $\mathcal{A}_{g,1}$ and $\mathcal{AB}_{g,1}$.

\subsection{Homology and homotopy actions}\label{subsec actions}

According to Griffith \cite{Gri}, the subgroup $\mathcal{B}_{g,1}$ (resp. $\mathcal{A}_{g,1}$) is characterized by the fact that its action on $\pi_1(\Sigma_g,x_0)$  preserves the normal subgroup generated by the curves $\beta_1, \dots, \beta_g$ (resp. $\alpha_1, \dots, \alpha_g$). As a consequence the action on homology of $\mathcal{B}_{g,1}$ (resp. $\mathcal{A}_{g,1}$ ) preserves the Lagrangian $B$, (resp. $A$). 

If one writes the matices of the symplectic group $\text{Sp} \omega$ as blocks according to the decomposition $H = A \oplus B$, then the image of $\mathcal{B}_{g,1} \rightarrow \text{Sp} \omega$ is contained in the subgroup $\text{Sp}_B \omega$ of matrices of the form~:
$\left( \begin{matrix}
G_1 & 0 \\ M & G_2 
        \end{matrix}
\right)$. 

Such matrices are symplectic if and only if $G_2 = ^tG_1^{-1}$ and $ ^tG_1 M$ is symmetric and we have an isomorphism~:
\[
\begin{array}{ccl}
\text{Sp}_B \omega & \stackrel{\sim}{\longrightarrow} & \text{GL}_g(\Z) \ltimes S_{g}(\Z) \\
\left( \begin{matrix}
G & 0 \\ M & ^tG^{-1} 
        \end{matrix}
\right) & \longmapsto & (G , ^tG M).
\end{array}
\]
Here $S_{g} (\Z)$ denotes the symmetric group of $g \times g$ matrices over the integers; the composition on the semi-direct product is given by the rule $(G, S) (H, T) = (GH, ^tH S H + T)$.
Checking on generators (see Suzuki \cite{Suz}) of $\mathcal{B}_{g,1}$ we get: 
\begin{lemma}\label{lem h1(bg,1))}
There is a short exact sequence of groups~:
\[
\xymatrix@1{ 1 \ar[r] & \mathcal{TB}_{g,1} \ar[r] & \mathcal{B}_{g,1}  \ar[r] & {\mathrm{GL}_g(\Z)} \ltimes S_{g}(\Z) \ar[r] & 1.}
\]
\end{lemma}

An analogous statement holds for $\mathcal{A}_{g,1}$ replacing the lagrangian $B$ by $A$.

We also recall a result due to Luft \cite[Corollary 2.1]{MR0500977}

\begin{lemma}\label{lem Luftonpi}
The natural homomorphism $\mathcal{B}_{g,1} \rightarrow  \text{Aut } \pi_1(\mathcal{H}_g, x_0)$ is onto.
\end{lemma}

Again an analogous statement holds for $\mathcal{A}_{g,1}$. If we restrict our attention to $\mathcal{AB}_{g,1}$ then the natural homomorphism $\mathcal{AB}_{g,1} \rightarrow \text{Aut } \pi_1(\mathcal{H}_g,x_0)$ is still an automorphism for the elements of $\mathcal{B}_{g,1}$ that hit the generators of the automorphism group in Luft's papper are readily seen to live in fact in $\mathcal{AB}_{g,1}$ (for a geometric description of these generators see \cite{MR0500977} or \cite{Suz} and for an algebraic description see Section \ref{sec generatosLT}). As for the previous Lemma, checking on generators we get~:

\begin{lemma}\label{lem actAB}
There is a short exact sequence of groups~:
\[
\xymatrix@1{ 1 \ar[r] & \mathcal{TAB}_{g,1} \ar[r] & \mathcal{AB}_{g,1} \ar[r] &  \mathrm{GL}_g(\Z) \ar[r] & 1.}
\]
\end{lemma}

\subsection{Heegaard splittings of homology spheres}\label{sub heegarrdsplit}

It is well known that by glueing two handlebodies with opposite orientations along a diffeomorphism of their boun\-da\-ry  one can construct all oriented compact $3$ manifolds. Choose a map $\iota_g \in \mathcal{M}_{g,1}$ such that $S^3 = \mathcal{H}_g \cup_{i_g} -\mathcal{H}_g$. If we twist this glueing by an arbitrary map in $\mathcal{T}_{g,1}$ we get back a new homology sphere $S^3_\phi = \mathcal{H}_g \cup_{i_g \phi} -\mathcal{H}_g$ and in fact we can get all homology spheres by letting $g$ vary \cite{Mo1}. More precisely,  consider the following equivalence relation on $\mathcal{T}_{g,1}$:
$$
\phi \sim \psi \Leftrightarrow \exists \zeta_a \in \mathcal{A}_{g,1} \exists \zeta_b \in \mathcal{B}_{g,1} \mbox{ such that } \zeta_b \phi \zeta_a = \psi.\quad (\ast)
$$

Moreover define the  \emph{stabilization map} on the mapping class group as follows. Glue one of the boundary components of a two holed torus on the boundary of $\Sigma_{g,1}$  to get $\Sigma_{g+1,1}$. Extending an element of $\mathcal{M}_{g,1}$ by the identity over the torus yields an injective homomorphism $\mathcal{M}_{g,1} \hookrightarrow \mathcal{M}_{g+1,1}$, this is the stabilization map. This map is compatible with the action on homology and is compatible with the definition of the above two subgroups $\mathcal{A}_{g,1}$ and $\mathcal{B}_{g,1}$. In particular the equivalence relation  $(\ast)$ is compatible with the stabilisation map . It is also possible to choose the map $i_g$ to be compatible with the stabilization map $i_{g+1} \vert_{\Sigma_{g,1}} = i_g$ and we have the following precise version of Heegaard splittings of integral homology spheres (see \cite{Mo1} for a proof):

\begin{theorem} 
The following map is well defined and is bijective:
$$
\begin{array}{ccc}
\displaystyle \lim_{g \to \infty} \mathcal{T}_{g,1}/\sim &\longrightarrow & \mathcal{S}(3)\\
\phi & \longmapsto & S^3_\phi = \mathcal{H}_g \cup_{ \phi} -\mathcal{H}_g
\end{array}
$$
\end{theorem}

From a group-theoretical point of view the equivalence relation $(\ast)$ is quite unsatisfactory, for it looks like, but is not, a double coset relation on the Torelli group. In fact it is the composite of a double coset relation in the Torelli group and a conjugacy-induced equivalence relation:

\begin{theorem}\label{thm nllerel}
Two maps $\phi, \psi \in \mathcal{T}_{g,1}$ are equivalent if
and only if there exists a map $\mu \in \mathcal{AB}_{g,1}$ and two maps $\xi_a \in \mathcal{TA}_{g,1}$ and
$\xi_b \in \mathcal{TB}_{g,1}$ such that $\phi = \mu \xi_b \psi \xi_a
\mu^{-1}$.
\end{theorem}

\begin{proof}
The ``if'' part of the theorem is trivial. Conversely, assume that $\psi = \xi_b \phi
\xi_a$, where $\psi, \phi \in \mathcal{T}_{g,1}$. Projecting this equality on
$\text{Sp} \omega$ we get $Id = \mathrm{H}_1 (\xi_a)  \mathrm{H}_1(\xi_b)$. According to Lemma \ref{lem h1(bg,1))} the matrix  $\mathrm{H}_1(\xi_b)$ is of the form  
\[
\left( \begin{matrix} G & 0 \\ M & ^tG^{-1} \end{matrix}
\right).
\]
Similarly, $\mathrm{H}_1(\xi_a)$ is of the form 
\[\left( \begin{matrix} H & N \\
0 &  ^tH^{-1} \end{matrix} \right).
\]
Therefore~:
\[
Id =  \left( \begin{matrix} H & N \\
0 &  ^tH^{-1} \end{matrix} \right) \left( \begin{matrix} G & 0 \\ M & ^tG^{-1} \end{matrix}
\right) =  \left( \begin{matrix}
HG + NM & N^tG^{-1} \\
^tH^{-1} M & ^tH^{-1} G^{-1} \end{matrix} \right).
\]
Thus~:
\[
N = 0 = M \text{ and } G= H^{-1}.
\]

In particular  \(\mathrm{H}_1(\xi_a),\ \mathrm{H}_1(\xi_b) \in \mathrm{GL}(g, \mathbf{Z}) 
\) and \( \mathrm{H}_1(\xi_a) = \mathrm{H}_1(\xi_b)^{-1} \). By Lemma \ref{lem actAB} we can choose
a map $\mu \in \mathcal{AB}_{g,1}$ such that $\mathrm{H}_1(\mu) = \mathrm{H}_1(\xi_b)$,
and we get
\[
\psi = \mu \circ (\mu^{-1} \xi_b) \phi (\xi_a \mu) \circ \mu^{-1}. 
\]
By construction $(\mu^{-1} \xi_b) \in \mathcal{TB}_{g,1}$
and $(\xi_a \mu) \in \mathcal{TA}_{g,1}$.
\end{proof}

\section{Trivial cocycles and invariants}\label{sec Trivcocycle}

\subsection{The Johnson homomorphism}\label{subsec johnsonhomom}

Computing the action of the Torelli group on the second nilpotent quotient of $\pi_1(\Sigma_{g,1},x_0)$ Johnson defines a morphism of groups known as the first Johnson homomorphism:
\[
\tau : \mathcal{T}_{g,1} \longrightarrow \wedge^3 H.
\]

Notice that the mapping class group $\mac$ acts naturally by conju\-ga\-tion on $\mathcal{T}_{g,1}$ and acts also on $\wedge^3 H$ via its natural action on homology. In \cite{Joh1}, \cite{Joh2},\cite{Joh3} Johnson proves that

\begin{proposition}\label{prop johnsonhom}
The map $\tau$ is $\mac$-equivariant with respect to the above actions. Up to finite dimensional $\mathbf{Z}/2 \mathbf{Z}$-vector space $\wedge^3 H$  is the abelianization of the Torelli group: any homomorphism $\mathcal{T}_{g,1} \rightarrow A$ where $A$ is an abelian group without $2$-torsion factors uniquely through $\tau$.
\end{proposition}

\subsection{From invariants to trivial cocycles}\label{subsec invtococy}

Consider an integer valuated invariant of homology spheres $F~: \mathcal{S}(3) \rightarrow \Z$. Precomposing with the canonical maps $\mathcal{T}_{g,1} \rightarrow \lim_{g \to \infty} \mathcal{T}_{g,1}/\sim \rightarrow \mathcal{S}(3)$ we get a family of maps $F_g : \mathcal{T}_{g,1} \rightarrow \Z$. Since the stabilization maps are injective  the map $F_g$ determines by restricton all maps $F_{g'}$ for $g'< g$. Therefore, as stated in the introduction, we avoid the peculiarites of the first Torelli groups by restricting ourselves to $g \geq 3$. We also consider the associated trivial cocycles, which measure the failure of the maps $F_g$ to be homomorphisms of groups
\[
\begin{array}{ccl}
C_g : \mathcal{T}_{g,1} \times \mathcal{T}_{g,1} & \longrightarrow & \Z \\
(\phi,\psi) & \longmapsto & F_g(\phi) + F_g(\psi) -F_g(\phi\psi).
\end{array}
\]

Since $F$ is an invariant the cocycles $C_g$ inherit the following pro\-per\-ties:

\begin{enumerate}
\item The cocycles $C_g$ are compatible. i.e. the following diagram of maps com\-mutes: 
\[
\xymatrix{
\mathcal{T}_{g,1} \times \mathcal{T}_{g,1} \ar@{^{(}->}[r] \ar[dr]_-{C_g} & \mathcal{T}_{g+1,1} \times \mathcal{T}_{g+1,1} \ar[d]^-{C_{g+1}} \\
 & \Z.
}
\]
\item The cocycles $C_g$ are invariant under conjugation by elements in $\mathcal{AB}_{g,1} : C_g(\phi - \phi^{-1}, \phi - \phi^{-1}) = C_g(-,-).$
\item If $\phi \in \mathcal{TB}_{g,1}$ or $\psi \in \mathcal{TA}_{g,1}$ then $C_g(\phi,\psi) = 0$.
\end{enumerate}

\begin{proposition}\label{prop cgnon0}
The cocycle $C_g$ is constantly equal to $0$ if and only if $F_g$ is the zero map.
\end{proposition}
\begin{proof}
If $C_g$ is $0$ then $F_g$ is a morphism of groups and therefore factors via $\tau$:
\[
\xymatrix{
\mathcal{T}_{g,1} \ar[d]^-{\tau} \ar[dr]^-{F_g} & \\
\wedge^3 H \ar[r]^-{\overline{F}_g} & \Z.
}
\] 
The morphism $\overline{F}_g$ is  then $\mathrm{GL}_g(\Z) = \mathcal{AB}_{g,1}/\mathcal{TAB}_{g,1}$-invariant. As $-Id \in \mathrm{GL}_g(\Z)$ acts as $-Id$ on $\wedge^3 H$, we get that $\overline{F}_g= 0$.
\end{proof}

Any two trivializations of a given trivial cocycle differ by a homomorphism of groups  and by the same argument we get~:

\begin{proposition}\label{prop unicite}
Any family of trivial cocycles satisfying $(1)-(3)$ coresponds to at most one invariant of homology spheres.
\end{proposition}

\subsection{From trivial cocycles to invariants}\label{subsec cocytoinv}

Conversely, what are the conditions for a family of trivial $2$-cocycles $C_g$ on $\mathcal{T}_{g,1}$ satisfying properties $(1)-(3)$ to actually provide an invariant ?

Firstly we need to check the existence of an $\mathcal{AB}_{g,1}$-invariant tri\-via\-lization of each $C_g$. This is a cohomological problem.

Denote by $\mathcal{Q}_{C_g}$ the set of all trivializations of the cocycle $C_g$: 
\[
\mathcal{Q}_{C_g} = \{ \xymatrix@1{q~:\mathcal{T}_{g,1} \ar[r] & \Z} \ | q(\phi) + q(\psi) - q(\phi \psi) = C_q(\phi,\psi)\}.
\]
Recall that any two trivializations differ by an element of the group $\mathrm{Hom}(\mathcal{T}_{g,1}, \Z) = \mathrm{Hom}(\wedge^3 H, \Z).$ As the cocycle $C_g$ is invariant under conjugation by $\mathcal{AB}_{g,1}$ this later group acts on $\mathcal{Q}_{C_g}$ via its conjugation action on the Torelli group.  Explicitly if $\phi \in \mathcal{AB}_{g,1}$ and $q \in \mathcal{Q}_{C_g}$ then $\phi \cdot q (\eta)  = q(\phi \eta \phi^{-1}).$ This action confers the set $\mathcal{Q}_{C_g}$ the structure of an affine set over the abelian group $\mathrm{Hom}(\wedge^3 H, \Z)$. Choose an arbitrary element $q \in  \mathcal{Q}_{C_g}$ and define a map as follows
\[
\begin{array}{rcl}
\rho_q~: \mathcal{AB}_{g,1} & \longrightarrow & \mathrm{Hom}(\wedge^3 H, \Z) \\
\phi & \longmapsto & \phi \cdot q - q.
\end{array}
\]
A direct computation shows that $\rho_q$ is a derivation and that the difference $\rho_q - \rho_{q'}$ for two elements in $\mathcal{Q}_{C_g}$ is a principal derivation, therefore we have a well defined cohomology class 
$$\rho(C_g) \in \mathrm{H}^1(\mathcal{AB}_{g,1}; \mathrm{Hom}(\wedge^3 H, \Z))$$
called the \emph{torsor} of the cocycle $C_g$.

By construction if the action of $\mathcal{AB}_{g,1}$ on $\mathcal{Q}_{C_g}$ has a fixed point the class $\rho(C_g)$ is trivial. Conversely, if $\rho(C_g)$ is trivial, then for any $q \in \mathcal{Q}_{C_g}$ the map $\rho_q$ is a principal derivation: there exists $m \in  \mathrm{Hom}(\mathcal{T}_{g,1}, \Z)$ such that
\[
\forall \phi \in \mathcal{AB}_{g,1} \ \rho_q(\phi) = \phi \cdot m - m.
\]
In particluar the element $q - m \in \mathcal{Q}_{C_g}$ is fixed under the action of $\mathcal{AB}_{g,1}$.

\begin{proposition}\label{prop torseurs}
The natural action of $\mathcal{AB}_{g,1}$ on $\mathcal{Q}_{g,1}$ admits a fixed point if and only if the associated torsor $\rho(C_g)$ is trivial.
\end{proposition}

Arguing as in Proposition \ref{prop cgnon0} one checks that the $\mathcal{AB}_{g,1}$-invariant trivialization of $C_g$, if it exists, is unique. If we have fixed points $q_g$ for all $g$, by unicity, we have that $q_{g+1}$ restricted to $\mathcal{T}_{g,1}$ is equal to $q_g$. Therefore we have a well-defined map
\[
\displaystyle  q =\lim_{g \to \infty} q_g~: \lim_{g \to \infty } \mathcal{T}_{g,1} \longrightarrow \Z. 
\]

This is the only candidate to be an invariant of homology spheres. For this map to be an invariant, since it is already $\mathcal{AB}_{g,1}$-invariant, we only have to prove that it is constant on the double cosets $\mathcal{TB}_{g,1} \backslash \mathcal{T}_{g,1} \slash \mathcal{TA}_{g,1}$.

From property $(3)$ of our cocycle we get that $\forall \phi \in  \mathcal{T}_{g,1}, \forall \psi_a \in \mathcal{TA}_{g,1}$ and $\forall \psi_b \in \mathcal{TB}_{g,1}$:
\begin{eqnarray*}
q_g(\phi) - q_g(\psi_b \phi) & = & -q_g(\psi_b) \\
q_g(\phi) - q_g(\phi \psi_a) & = & -q_g(\psi_a) \\
\end{eqnarray*}

\begin{theorem}\label{thm nulliteqg}
For each $g \geq 3$ the induced homomorphisms 
\[ q_g~: \mathcal{TB}_{g,1} \rightarrow \Z \text{ and } q_g~: \mathcal{TA}_{g,1} \rightarrow \Z \] are trivial.
\end{theorem}

\begin{proof} We only give the proof for the morhism $q_g~: \mathcal{TB}_{g,1} \rightarrow \Z$, the other case is similar. 

Denote by $ \mathcal{L}_{g,1}$ the kernel of  the map $\mathcal{B}_{g,1} \twoheadrightarrow \text{Aut } \pi_1(\mathcal{H}_g)$. This was identified by Luft \cite{MR0500977} as the ``Twist group" of the handlebody $\mathcal{H}_g$.

Consider the following commutative diagram
\[
\xymatrix{
1 \ar[r] & \mathcal{L} \cap \mathcal{TB}_{g,1} \ar[r] \ar@{^{(}->}[d] & \mathcal{TB}_{g,1} \ar[r] \ar@{^{(}->}[d] & IA \ar[r] \ar@{^{(}->}[d] & 1\\
1 \ar[r] & \mathcal{L}_{g,1} \ar[r] & \mathcal{B}_{g,1} \ar[r] & \text{Aut} \pi_1(\mathcal{H}_g) \ar[r] & 1
}
\]

Here $IA$ stands for the kernel of the natural map $\text{Aut} (\mathcal{H}_g) \rightarrow \mathrm{GL}_g(\Z).$ In Section \ref{sec generatosLT} we prove the following result:

\begin{proposition}\label{prop gendeLTB}
The group $\mathcal{L} \cap \mathcal{TB}_{g,1}$ is generated by maps of the form $T_\beta T_{\beta'}^{-1}$, where $\beta$ and $\beta'$ are two homologous non-isotopic and disjoint simple closed curves on $\Sigma_{g,1}$ such that each one bounds a properly embedded disc in $\mathcal{H}_g$.
\end{proposition}

\begin{corollary}\label{cor nulsurLTB}
The morphism $q_g : \mathcal{L} \cap \mathcal{TB}_{g,1} \rightarrow \Z$ is trivial.
\end{corollary}
\begin{proof}
By Proposition \ref{prop gendeLTB} it is enough to prove that $q_g$ vanishes on the aforementioned maps $T_\beta T_{\beta'}^{-1}$. As our embedding of $\Sigma_g$ into $\mathbf{S}^3$ is standard, there exists a simple closed curve, $\alpha \subset \Sigma_{g,1}$ which bounds a properly embedded disc on $-\mathcal{H}_g$ (the outer handlebody) and which intersects each of the curves $\beta$ and $\beta'$ in exactly one point. Consider a regular neighbourhood of the union $\beta \cup  \alpha \cup \beta'$, it is a $3$-ball, whose intersection with the surface looks like in Figure \ref{fig intersection}.

\begin{figure}
\input{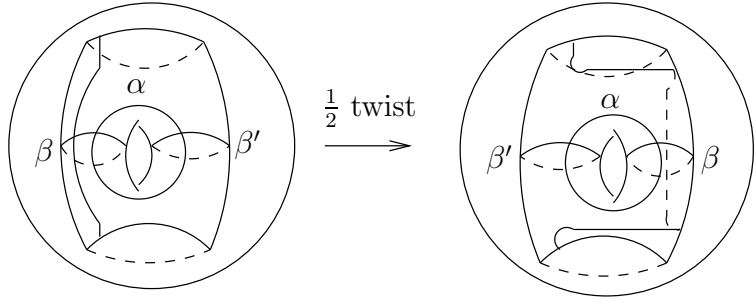}
\caption{Regular neighbourhood of $\beta \cup \alpha \cup \beta'$ and half twist}\label{fig intersection}
\end{figure}

There is a half twist map $\psi$ inside this ball that exchanges the curves $\beta$ and $\beta'$. This half twist map $\psi$ belongs to $\mathcal{AB}_{g,1}$ since it is a self diffeomorphism of the depicted $3$-ball and can be extended by the  identity outside this ball. In particular, since $q_g$ is $\mathcal{AB}_{g,1}$-invariant:
\begin{eqnarray*}
q_g(T_\beta T_{\beta'}^{-1}) & = & q_g (\psi T_\beta T_{\beta'}^{-1} \psi^{-1}) \\
& = & q_g (T_{\psi(\beta)} T_{\psi(\beta')}^{-1}) \\
& = & q_g(T_{\beta'} T_{\beta}^{-1}) \\
& = & - q_g(T_\beta T_{\beta'}^{-1})
\end{eqnarray*}
and therefore $q_g\vert_{\mathcal{L} \cap \mathcal{TB}_{g,1}} =0$.
\end{proof}

As a consequence $q_g$ factors through $IA$. As the action on the fundamental group of the inner handlebody $\mathcal{H}_g$ induces a surjective map $\mathcal{AB}_{g,1} \rightarrow \text{Aut } \pi_1 (\mathcal{H}_g)$ we can even view $q_g$ as an $\text{Aut } \pi_1 (\mathcal{H}_g)$-invariant map $q_g : IA \rightarrow \Z$. Let $\alpha_1, \dots, \alpha_g$ denote the generators of $\pi_1(\mathcal{H}_g)$ (this can be identified with the curves  in Figure \ref{fig model}). According to Magnus \cite{mag1934}, the group $IA$ is normally generated as a subgroup of $\text{Aut } \pi_1(\mathcal{H}_g)$ by the automorphism $K_{12}$ given by $K_{12}(\alpha_1) = \alpha_2 \alpha_1 \alpha_2^{-1}$ and $K_{12}(\alpha_i) = \alpha_i \text{ for } i \geq 2$. By invariance, $q_g$ is determined by its value on $K_{12}$.
An easy computation shows that if we denote by $\sigma_2$ the automorphism given by $\sigma_{2} (\alpha_2) = \alpha_2^{-1}$ and $\sigma_2(\alpha_i) = \alpha_i \text{ for } i \neq 2$  then $\sigma_2 K_{12} \sigma_2^{-1} = K_{12}^{-1}$. In particular $q_g(K_{12}) = q_g(\sigma_2 K_{12} \sigma_2^{-1}) = - q_g(K_{12})$ and $q_g$ is must be trivial.
\end{proof}

We sumarize the discussion of this section in the following

\begin{theorem}\label{thm main}
A family of cocycles $(C_g)_{g \geq 3}$ on the Torelli groups $\mathcal{T}_{g,1}$, $ g \geq 3$ satisfying conditions $(1)-(3)$ provides a compatible family of trivializations $F_g : \mathcal{T}_{g,1} \rightarrow \Z$ that reassemble into an invariant of homology spheres
\[
\displaystyle \lim_{g \to \infty} F_g : \mathcal{S}(3) \rightarrow \Z
\]
if and only if the following two conditions hold:
\begin{enumerate}
\item[(i)] The associated cohomology classes $[C_g] \in \mathrm{H}^2(\mathcal{T}_{g,1}; \Z)$ are trivial.
\item[(ii)] The associated torsors $\rho(C_g) \in \mathrm{H}^1(\mathcal{AB}_{g,1}, \wedge^3 \mathrm{H}_1(\Sigma_g))$ are trivial.
\end{enumerate}
In this case the maps $F_g$ are the unique $\mathcal{AB}_{g,1}$-invariant trivializations of the cocycles $C_g$.
\end{theorem}

\section{Application to the Casson invariant}\label{sec Casson}

If one is interested in invariants that come from pull-backing cocycles defined on abelian groups without $2$-torsion, in view of Propositon \ref{prop johnsonhom} it is enough to study the case where the abelian group is $\wedge^3 H$ and the homomorphism is $\tau$.

Recall that we have a decomposition $H = A \oplus B$, this induces the decomposition $\wedge^3 H = \wedge^3 A \oplus B \wedge \wedge^2 A + A \wedge \wedge^2 B \oplus \wedge^3 B$.  Set $W_A = \wedge^3 A$, $W_B = \wedge^3 B$ and $W_{AB} = \oplus B \wedge \wedge^2 A + A \wedge \wedge^2 B$. The Johnson homomorphism computes the action of the Torelli group on the second nilpotent quotient of the fundamental group of $\Sigma_{g,1}$. Computing on specific eleemnts one can chek that (see \cite{Mo1}) :

\begin{lemma}\label{lem imTbtau1}
The image of $\mathcal{TB}_{g,1}$ under $\tau$ in $\wedge^3 H$ is $W_A \oplus W_{AB}$, the image of $\mathcal{TA}_{g,1}$ is $W_{AB} \oplus W_B$. 
\end{lemma}

For each $g$, the symplectic pairing $\omega : \wedge^2 H \rightarrow \Z$ induces  natural pairings  $J_g~: W_A \wedge W_B \rightarrow \Z$ and $^t J_g~: W_B \wedge W_A \rightarrow \Z$  that can be naturally extended to  degenerate bilinear forms on $\wedge^3 H= W_A \oplus W_{AB} \oplus W_B$, with matrices:
\[
J_g := \left(\begin{matrix}
0 & 0 & 0 \\
0 & 0 & 0 \\
Id & 0 & 0
     \end{matrix}
\right), \
^tJ_g := \left(\begin{matrix}
0 & 0 & Id \\
0 & 0 & 0 \\
0 & 0 & 0
     \end{matrix}
\right).
\]
Notice that bilinear forms are naturally $2$-cocycles on abelian groups. 

\begin{proposition}\label{prop unicitecocycleab}
For each $g \geq 3$, the cocycle $J_g$ is the unique cocycle (up to a multiplicative constant) on $\wedge^3 H$ which pull-back on the Torelli group $\mathcal{T}_{g,1}$ satisfies conditions $(2)$ and $(3)$. Moreover once we have fixed a common multiplicative constant the family of pull-backed cocycles satisfies also $(1)$.
\end{proposition}
\begin{proof}
Fix an integer $ n \in \mathbf{Z}$. It is obvious from the definition and from Lemma \ref{lem imTbtau1} that the family $(nJ_g)$ satisfies $(1)$, $(2)$ and $(3).$ 

Let $B$ denote an arbitrary cocycle on $\wedge^3 H$ which pull-back on $\mathcal{T}_{g,1}$ satisfies $(2)$ and $(3)$.

Write each element $w  \in  \wedge^3 H$ as $w_a + w_{ab} + w_b$ according to the decomposition $W_A \oplus W_{AB} \oplus W_B$. The cocycle relation together with condition $(3)$ and Lemma \ref{lem imTbtau1} imply that
\[
\forall v,w \in \wedge^3 H \quad B(v,w) = B(v_a,w_b).
\]

We first prove that $B$ is bilinear. For the linearity on the first variable compute
\begin{eqnarray*}
B(u+v,w) & = & B(u_a + v_a, w_b) \\
 &= & B(v_a, w_b) + B(u_a, v_a + w_b) - B(u_a, v_a) \\
 & = & B(u_a,w_b) + B(v_a, w_b) \\
  & = & B(u,w) + B(v,w).
\end{eqnarray*}
A similar proof holds for the linearity on the second variable. 

By the equivariance properties of $\tau$, the subgroup $\mathcal{AB}_{g,1} \subset \mac$ acts on $\wedge^3 H$ via the projection $\mathcal{AB}_{g,1} \stackrel{\mathrm{H}_1}{\longrightarrow} \mathrm{GL}_g(\Z)$. It is well known that the only $\mathrm{GL}_g(\Z)$-invariant bilinear forms on $\wedge^3 (A \oplus B)$ are the pairing $J_g$ and the dual pairing $^t J_g$, so that $B = n J + m {}^tJ$ for some integers $n$ and $m$. As condition $(3)$ implies that $\tau(\mathcal{TB}_{g,1}) = W_B \oplus W_{AB}$ has to be in the kernel of $B$, evaluating on the elements of $W_B$ yields that $m =0$.
\end{proof}

We would like to apply Theorem \ref{thm main} to the family $(J_g)$ or to one of its multiples. First we must check the triviality of the cocycles.

\begin{proposition}\label{prop triv2J}
For all $g \geq 3$, the pull-back of the $2$-cocycle $2J_g$ is trivial.
\end{proposition}
\begin{proof}
From Johnson(see Proposition \ref{prop johnsonhom}) we know that the abe\-lia\-nization of $\mathcal{T}_{g,1}$ is equal to $\Lambda^3 H \oplus V$ where $V$ is a $\Z /2\Z$-vector space. By the universal coefficients theorem $\HH^2(\mathcal{T}_{g,1}, \Z) = \mathrm{Hom}(\HH_2(\mathcal{T}_{g,1},\Z),\Z) \oplus \text{Ext}^1(\HH_1(\mathcal{T}_{g,1},\Z), \Z)$. The first factor is torsion free and the second factor is isomorphic to $\text{Ext}^1(V,\Z)$ which is a $\Z/ 2\Z$ vector space .

By naturality we have a commutative diagram for cohomology groups with trivial coefficients~:
\[
\xymatrix{
{\HH^2(\wedge^3 H~; \mathbf{Q})} \ar^{\tau_{ \mathbf{Q}}^\ast}[r] & {\HH^2(\mathcal{T}_{g,1}~; \mathbf{Q})} \\
{\HH^2(\wedge^3 H~; \mathbf{Z})} \ar^{\tau^\ast}[r] \ar@{_{(}->}[u]& {\HH^2(\mathcal{T}_{g,1}~; \mathbf{\mathbf{Z}})} \ar[u].
}
\]

R. Hain \cite{Hain} has computed the kernel of $\tau_{\mathbf{Q}}^\ast$ and a direct computation shows that the image of the pull-back of $J_g$ in $\HH^2(\wedge^3 H; \mathbf{Q})$ lies in this kernel. In particular $\tau^\ast(J_g)$ is anihilated by multiplication by $2$ so that the pull-back of the $2$-cocycle $2J_g$ is trivial. 
\end{proof}

We will come back to the homology class of the pull-back of the cocycle $J_g$ later (see Proposition \ref{prop classecasson}). To avoid an unnecessarily heavy notation from now on we will also denote by $J_g$ the pull-back of the cocyle $J_g$ along the morphism $\tau$.  To see if there is an invariant associated to the family $(2J_g)_{g \geq 3}$ we have to check the triviality of the associated torsors:

\begin{proposition}\label{prop trivialtorsors}
For each $g \geq 3$, the torsor \[\rho(2J_g) \in \mathrm{H}^1(\mathcal{AB}_{g,1}; \mathrm{Hom}(\Lambda^3 H, \mathbf{Z}))\] is trivial.
\end{proposition}
\begin{proof}
By definition we have an exact sequence 
\[
\xymatrix{
1 \ar[r] & \mathcal{TAB}_{g,1} \ar[r] & \mathcal{AB}_{g,1} \ar[r] & \mathrm{GL}_g(\mathbf{Z})\ar[r] & 1
}
\]
We get an induced exact sequence in low-dimensional cohomology~:
\[
\xymatrix@C=.45cm{
0 \ar[r] & \mathrm{H}^1(\mathrm{GL}_g(\mathbf{Z});\mathrm{Hom}(\Lambda^3 H,\mathbf{Z})) \ar[r] & \mathrm{H}^1(\mathcal{AB}_{g,1}; \mathrm{Hom}(\Lambda^3 H, \mathbf{Z})) \ar[d] \\
& & \mathrm{H}^1(\mathcal{TAB}_{g,1};\mathrm{Hom}(\Lambda^3 H,\mathbf{Z}))^{\mathrm{GL}_g(\mathbf{Z})}.
}
\]

First we show that $\mathrm{H}^1(\mathrm{GL}_g(\mathbf{Z});\mathrm{Hom}(\Lambda^3 H,\mathbf{Z})) = 0$.

Let $f~: \mathrm{GL}_g(\mathbf{Z}) \longrightarrow \Lambda^3 H$ be any crossed morphism. As $-Id \in \mathrm{GL}_g(\mathbf{Z})$ acts as $-Id$ on  $\mathrm{Hom}(\Lambda^3 H,\mathbf{Z})$ and is central, for all $ S \in \mathrm{GL}_g(\mathbf{Z})$ we have $f(-Id \circ S) = f(-Id) - f(S) =  f(S) + S\cdot f(-Id)$. In particular, $\forall S \in \mathrm{GL}_g(\mathbf{Z})$, $2f(S) = f(-Id) - S\cdot f(-Id)$. Using the standard generators (elementary matrices) $E_{ij}$, defined by $E_{ij}(a_k) = a_k + \delta_{jk}a_i$ one shows that $f(-Id)$ is divisible by $2$, so $f$ itself is a principal derivation.

We are left with the exact sequence~:
$$
0 \rightarrow \mathrm{H}^1(\mathcal{AB}_{g,1}; \mathrm{Hom}(\Lambda^3 H, \mathbf{Z})) \stackrel{i_\ast}{\longrightarrow} \mathrm{H}^1(\mathcal{TAB}_{g,1};\mathrm{Hom}(\Lambda^3 H,\mathbf{Z}))^{\mathrm{GL}(g,\mathbf{Z})}.
$$
As $\mathcal{T}_{g,1}$ and therefore $\mathcal{TAB}_{g,1}$, acts trivially on the free abelian group $\mathrm{Hom}(\Lambda^3 H, \mathbf{Z})$, the universal coefficients theorem and the classical adjun\-ction properties of $\mathrm{Hom}$-groups gives us a canonical  $\mathrm{GL}(g, \mathbf{Z})$-equi\-va\-ri\-ant isomorphism
\[
\mathrm{H}^1(\mathcal{TAB}_{g,1};\mathrm{Hom}(\Lambda^3 H,\mathbf{Z})) \simeq  \mathrm{Hom}(\HH_1 \mathcal{TAB}_{g,1} \otimes \Lambda^3 H,\mathbf{Z}) 
\]

By construction of the torsor class, the image of $\rho(2J_g)$ under this isomorphism may be described as follows. Fix an arbitrary coboundary $q \in \mathcal{Q}_{2J_g}$. For each tensor $f \otimes l \in \HH_1(\mathcal{TAB}_{g,1}) \otimes \Lambda^3 H$, choose arbitrary lifts of $\phi \in \mathcal{TAB}(1)_{g,1}$ and  $\lambda \in \mathcal{T}_{g,1}$, then $\rho(2J_g)(f \otimes l) = q(\phi \lambda \phi^{-1}) - q (\lambda)$.
As $q$ is a coboundary of  $2J_g$ we get:
\begin{eqnarray*}
\rho(2J_g)(f \otimes l) 	& = & q(\phi \lambda \phi^{-1}) - q (\lambda) \\
			& = 	& q(\phi \lambda \phi^{-1} \lambda^{-1}) \\
			& =	& 2J_g(\tau(\phi), \tau(\lambda)) - 2J_g(\tau(\lambda), \tau(\phi)) \\
			& = 	& 0 \text{ by condition } ()2).
\end{eqnarray*}

\end{proof}

Applying Theorem \ref{thm main} we get
\begin{theorem}\label{thm construction}
The $\mathcal{AB}_{g,1}$-invariant trivializations of the pull-backs of the cocycles  $2J_g$ reassemble into an invariant of homology spheres $F~: \mathcal{S}(3) \rightarrow \Z$. Up to a multiplicative constant this is trivialization of a $2$-cocycle defined on an abelian group without $2$-torsion.
\end{theorem}

In the next paragraph   we identify the above invariant with the Casson invariant:

\begin{corollary}\label{cor caraccasson}
The Casson invariant is the unique integral valuated invariant of oriented homology $3$-spheres that comes from the tri\-via\-li\-zation of a $2$-cocycle defined on an abelian group without $2$-torsion.
\end{corollary}

We also get back one of the main results in \cite{Mo1}:

\begin{corollary}\label{cor cassonmorphisme}
Denote by $\mathcal{K}_{g,1}$ the kernel of the Johnson homomorphism $\tau$ (see \cite{Joh2} for geometric properties of this group). Then the Casson invariant restricted to $\mathcal{K}_{g,1}$ is a homomorphism of groups.
\end{corollary}

Assuming that we know that our invariant $F$ constructed out of the cocycles $2J_g$ is the opposite of the Casson invariant we proceed by showing that one can not get rid of the factor $2$. Denote the Casson invariant by $\lambda$ and by $\lambda_g$ if we view it as a function on $\mathcal{T}_{g,1}$.

\begin{proposition}\label{prop classecasson}
The pullback of the cocycle $J_g$ on the Torelli group defines a non-trivial cohomology class $[J_g] \in \mathrm{H}^2(\mathcal{T}_{g,1}; \Z)$ of order two. Moreover this classes are stable in the sense that the image of the class $[J_{g+1}]$ under the stabilisation map $\mathrm{H}^2(\mathcal{T}_{g+1,1}; \Z) \rightarrow \mathrm{H}^2(\mathcal{T}_{g,1}; \Z)$ is $[J_g]$.
\end{proposition}
\begin{proof}
If the pullbacks were trivial then the proof of Proposition \ref{prop trivialtorsors}
would carry on and provide us with an invariant $F_g ~: \mathcal{T}_{g,1} \rightarrow \Z$ associated to the family $(J_g)$. Then by the unicity of invariants associated to cocycles, Proposition \ref{prop unicite}, the invariant $2F_g$ would be the invariant associated to $2J_g$ so  we would have $2F_g = -\lambda_g$. Now, the Poincar\'e sphere has a Heegaard splitting of genus $2$ and therefore by stabilization, it has a Heegaard splitting of every genus $g \geq 3$. The Casson invariant of the Poincar\'e sphere is $1$ and therefore all functions $\lambda_g$ take the value $1$ and thus are not divisible by $2$.
\end{proof}

It is known that the mod $2$ reduction of the Casson invariant is the Rohlin invariant, which might be viewed as an homomorphism $R_g : \mathcal{T}_{g,1} \rightarrow \Z / 2\Z$ or equivalently as a cohomology  class $R_g \in \mathrm{H}^1(\mathcal{T}_{g,1}; \Z /2\Z) \simeq \mathrm{Hom}(\mathcal{T}_{g,1}, \Z /2\Z)$. By definition of the Bockstein homomorphism $\beta_\Z$ associated to the exact sequence \raisebox{.5ex}{\xymatrix@1{1 \ar[r] & \Z \ar[r]^-{\times 2} & \Z \ar[r] & \Z/ 2\Z \ar[r] &1 }}, we have~:

\begin{proposition}\label{prop bockstein}
For $g \geq 3$, the image of the class $R_g$ under the integral Bockstein $\beta_\Z : \mathrm{H}^1(\mathcal{T}_{g,1}; \Z/ 2\Z) \rightarrow \mathrm{H}^2(\mathcal{T}_{g,1}; \Z)$ is the non-trivial class $[J_g]$.
\end{proposition}
\subsection*{Identification of the Casson invariant}\label{subsec cassonprop}

In this section we prove that the inva\-riant $F$ derived from the familly of cocycles $(2J_g)$ is the Casson invariant $\lambda$. For the classical construction of this invariant and the study of some of its properties we refer the reader to the sruvey book \cite{MR1030042} or to the survey article \cite{MR1189008}. 

Additivity of the Casson invariant under connected sum is a consequence of Propositions \ref{prop surgery} and  \ref{prop boundarylink} hereafter (see \cite{MR1189008}). We give here an independent proof as it serves as an illustration of how to extract information on an invariant out of the associated cocycle. The statement can be easely reformulated for general invariants.

\begin{proposition}\label{prop connectedsum}
Denote by $M \sharp N$ the connected sum of the oriented manifolds $M$ and $N$. The invariant $F$ is additive with respect to connected sums of homology spheres: if $\phi \in \mathcal{T}_{g,1}$ and $\psi \in \mathcal{T}_{h,1}$ then
\[
F(\mathbf{S}^3_\phi \sharp \mathbf{S}^3_\psi) = F(\mathbf{S}^3_\phi) + F(\mathbf{S}^3_\psi).
\]
\end{proposition}
\begin{proof}
It is a classical result in the theory of Heegaard splittings that a Heegaard splitting of genus $g+ h$ for the connected sum $M \sharp N$ can be obtained from Heegaard splittings of  genus $g$ and $h$ respectively  of the manifolds $M$ and $N$ by splicing  along $3$-holed sphere, see Figure \ref{fig Heegsplit}.
\begin{figure}
\input{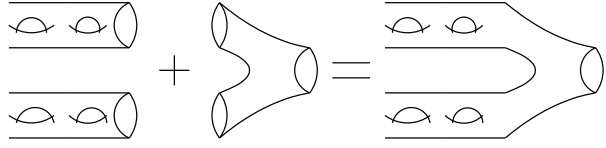}
\caption{Heegaard splitting of a connected sum}\label{fig Heegsplit}
\end{figure}

We see the surface $\Sigma_{g,1} \subset \Sigma_{g+h,1}$ as containing the $g$ first handles and the surface $\Sigma_{h,1} \subset \Sigma_{g+h,1}$ as containig the last $h$ handles. As for the stabilization map, extending by the identity we get two injective  maps $j_g~:\mathcal{T}_{g,1} \rightarrow \mathcal{T}_{g+h,1}$ and $j_h~:\mathcal{T}_{h,1} \rightarrow \mathcal{T}_{g+h,1}.$ If $\phi \in \mathcal{T}_{g,1}$ and $\psi \in \mathcal{T}_{h,1}$ then $\mathbf{S}^3_\phi \sharp \mathbf{S}^3_\psi = \mathbf{S}^3_{j_g(\phi) j_h(\psi)}.$
In particular
\[
F(\mathbf{S}^3_\phi \sharp \mathbf{S}^3_\psi) - F(\mathbf{S}^3_\phi) - F(\mathbf{S}^3_\psi) = -2J_{g+h}(\tau(j_g(\phi)), \tau(j_h(\psi))).
\]
If one turns back to the definition of the Johnson homomorphism \cite{Joh1} one can chek that $\tau(j_g(\phi)) \in \wedge^3 H$ only involves exterior powers of the homology classes $a_i, b_i$ for $i \leq g$ and $\tau(j_g(\psi)) \in \wedge^3 H$ involves only exterior powers of the homology classes $a_i, b_i$ for $g < i$. Therefore $J_{g+h}(\tau(j_g(\phi)), \tau(j_h(\psi))) = 0$.
\end{proof}

\begin{proposition}\label{prop surgery}
Let $M$ be an homology sphere and $K$ a knot in $M$. For an integer $n \geq 1$ denote by $K_{n}$ the result of performing a $\frac{1}{n}$-Dehn surgery on $K$. Then $F(K_{n+1}) -F(K_{n})$ is independent of $n$. 
\end{proposition}
\begin{proof}
There exists a Heegaard splitting for $M$ such that $K$ belongs to the surface $\Sigma_g$ and moreover is separating (i.e. $\Sigma_g - K$ has two connected components) (see \cite[Lemma 1.1 p.82 ]{MR1030042} for instance). Doing a $\frac{1}{n}$-Dehn surgery on $K$ is equivalent to modify the map $\phi \in \mathcal{T}_{g,1}$ such that $\mathbf{S}^3_\phi = M$ by the $n^{\text{th}}$ power of the Dehn twist along $K$.
Therefore 
\begin{eqnarray*}
F(K_{n+1}) - F(K_{n}) -F(\mathbf{S}^3_{T_K}) & = & F_{g}(T_K^{n+1}\phi) - F_g(T_K^{n} \phi)- F_g(T_K) \\
& = & -2J_g(\tau(T_K^{n} \phi),\tau(T_K)). 
\end{eqnarray*}

As $K$ is separating $\tau(T_K) = 0$ (see \cite{Joh1}), therefore $F(K_{n+1}) - F(K_{n})  = F(\mathbf{S}^3_{T_K})$ which is independent of $n$.
\end{proof}

Denote the quantity $F(K_{n+1}) - F(K_{n})$ by $F'(K)$ or $F'(K \subset M)$ if there is an ambiguity on the ambient space. If $(K,L)$ is a link in the homology sphere $M$, Proposition \ref{prop surgery} shows that $F(K_{k+1},L_{l+1}) - F(K_k,L_{l+1}) - F(K_{k+1},L_l) + F(K_k,L_l)$ is equal to $F'(K \subset L_{l+1}) - F'(K \subset L_l)$ and to $F'(L \subset K_{k+1}) - F'(L \subset K_l)$ and is therefore independent from $k$ and $l$. We denote this number by $F''(K,L)$.

\begin{proposition}\label{prop boundarylink}
Let $M$ be an homology sphere and let $K$ and $L$ be two knots in $M$ that bound disjoint Seifert surfaces in $M$. Then $F''(K,L) =0$.
\end{proposition}
\begin{proof}
There is a Heegaard splitting of $M$, say of genus $g$,  such that both $K$ and $L$ are separating curves on $\Sigma_g$( see \cite[Proposition 6.1 p.126]{MR1030042}). Denote by $\phi$ a map in $\mathcal{T}_{g,1}$ such that for the chosen Heegaard splitting $M = \mathbf{S}^3_\phi$. Recall that $T_K$ and $T_L$ belong to the kernel of the map $\tau : \mathcal{T}_{g,1} \rightarrow \wedge^3 H$ and that $F_g(\phi) + F_g(\psi) - F_g(\phi\psi) = 2J_g(\tau(\phi),\tau(\psi)).$
\begin{eqnarray*}
& & F(K_{k+1}L_{l+1}) - F(K_k L_{l+1}) - F(K_{k+1}L_l) + F(K_k L_l) \\
& = & F_g(T_L^{l+1}T_K^{k+1} \phi) - F_g(T_L^{l+1}T_K^{k} \phi) - F_g(T_L^{l}T_K^{k+1} \phi) + F_g(T_L^{l}T_K^{k} \phi) \\
& = & F_g(T_L^{l+1}T_K^{k+1}) + F(\phi) - F_g(T_L^{l+1}T_K^{k}) - F(\phi) - F_g(T_L^{l}T_K^{k+1}) \\
& &   - F(\phi) + F_g(T_L^{l}T_K^{k}) + F(\phi) \\
& = & F_g(T_L^{l+1}T_K^{k+1})- F_g(T_L^{l+1}T_K^{k}) - F_g(T_L^{l}T_K^{k+1})
 + F_g(T_L^{l}T_K^{k}) \\
& = & F_g(T_K) - F_g(T_L^{l}T_K^{k+1}) + F_g(T_L^{l}T_K^{k}) \\
& = & 2J_g(\tau(T_L^l), \tau(T_K)) \\
& = & 0.
\end{eqnarray*}
 
\end{proof}

Now from \cite[Proposition 1.3 p.235]{MR1189008} we learn  that up to a multiplicative constant there exists a unique invariant of oriented homology sphere that is $0$ on $\mathbf{S}^3$ and that satisfies Propositions \ref{prop surgery} and \ref{prop boundarylink} above. In particular as a consequence of the above two Propositions we get (see \cite[Proposition 1.3 p.235]{MR1189008})
\begin{enumerate}
\item (Surgery formula) Denote by $T$ the trefoil knot in $\mathbf{S}^3$ and by  $\Delta_K(t)$ the Alexander polynomial of a knot $K$ normalized at $1$. 

Then $F'(K) = \frac{1}{2}\Delta''_K(1)F'(T)$.
\item (Change of orientation) If $-M$ denotes the oriented homology sphere $M$ with opposite orientation, then $F(-M) = -F(M)$.
\end{enumerate}

\begin{corollary}\label{cor casson}
Up to a non-zero multiplicative constant, the invariant $F_g$ constructed in Theorem \ref{thm construction} is equal to the Casson invariant.
\end{corollary}

To fix the multiplicative constant we have to evaluate our invariant on one homology sphere and compute the value of the Casson invariant on the same sphere. Our invariant $F$ is most easely computed on maps of the form $\phi_a \phi_b$, with $\phi_a \in \mathcal{TA}_{g,1}$ and $\phi_b \in \mathcal{TB}_{g,1}$, for then:
\begin{eqnarray*}
F(\phi_a \phi_b) & = & -2J(\tau(\phi_a), \tau(\phi_b)) + F(\phi_a) + F(\phi_b) \\
& = & -2J(\tau(\phi_a), \tau(\phi_b)).
\end{eqnarray*}

\begin{figure}
\input{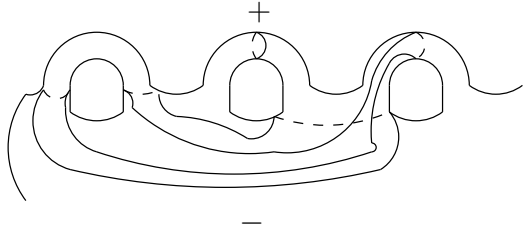}
\caption{Supporting curves for the map $\phi_b$}\label{fig phib}
\end{figure}

Consider the Dehn twists around the curves in Figures \ref{fig phia} and \ref{fig phib} where $+$ means a Dehn twist and $-$ the inverse of a Dehn twist.
\begin{figure}
\input{phia.pstex_t}
\caption{Supporting curves for the map $\phi_a$}\label{fig phia}
\end{figure}

The homology of the curves in Figure \ref{fig phia} is $\pm a_2$, and the spine of the genus one sub-surface they bound has homology $a_1$ and $b_1 - a_3$. Therefore, according to Johnson \cite{Joh1} $\tau(\phi_a)) = a_1 \wedge(b_1 - a_3) \wedge a_2$.
Similarly, the curves in figure \ref{fig phib} have homology $\pm b_2$, the spin of the subsurface of genus one they bound has homology $a1 + b_3, b_1$, so that $\tau(\phi_b)) = (a_1 +b_3) \wedge b_1 \wedge b_2$.

From the definition of the cocycle $J_g$ we get that 
\[
F(S^3_{\phi_a \phi_b}) = -2J_g(a_1 \wedge(b_1 - a_3) \wedge a_2,(a_1 +b_3) \wedge b_1 \wedge b_2 ) = 2.
\]

To compute the Casson invariant of $\mathbf{S}^3_{\phi_a\phi_b}$ we proceed as follows. We have expressed the map $\phi_a \phi_b$ as a product of Dehn twists $T_{\gamma_4} T_{\gamma_3}^{-1} T_{\gamma_2} T_{\gamma_1}^{-1}$. Choose a small $\varepsilon> 0$ and push the first curve $\gamma_1$ by epsilon in the inner handlebody, similarly push the curve $\gamma_2$ by $2\varepsilon$ and so on. This yelds a link in the inner handlebody with each component marked by $\pm 1$ according to the corresponding Dehn twist. Viewing this link in $\mathbf{S}^3$ we get a surgeery description of $\mathbf{S}^3_{\phi_a\phi_b}$ (see for instance \cite[p. 275]{MR1277811}). A straightforward computation using the surgery formula shows that the Casson invariant of the homology sphere $S^3_{\phi_a \phi_b}$ is $-2$.

\begin{proposition}\label{propo moritadef}
Denote by $\lambda~:\mathcal{S}(3) \rightarrow \Z$ the Casson invariant. Then the functions $\lambda_g$ for $g \geq 3$ satisfy the equation:
\[
\forall \phi, \psi \in \mathcal{T}_{g,1},  \ \ \lambda(\phi \psi) - \lambda_g(\phi) - \lambda_g(\psi) = 2J_g(\tau(\phi), \tau(\psi)).
\]
\end{proposition}

\begin{corollary}\label{cor identification}
The invariant $F$ constructed in Theorem \ref{thm construction} is the Casson invariant.
\end{corollary}

\section{Generators for the Luft-Torelli group}\label{sec generatosLT}

In this section we finally prove Proposition \ref{prop gendeLTB}. Before we need to recall some known facts on Dehn twists and maps in $\mathcal{TB}_{g,1}$.

In \cite{MR0500977} Luft identified the kernel $\mathcal{L}_{g,1}$ of the map \[ \xymatrix@1{ \mathcal{B}_{g,1} \ar[r] & \text{Aut } \pi_1(\mathcal{H}_g) \ar[r] & 1}\] with the so-called ``Twist group": the subgroup of $\mac$ generated by Dehn twists around simple closed curves that are contractible in $\mathcal{H}_{g}$.

In analogy with the generators of the Torelli group defined by Johnson \cite{Joh3}, we define a Contractible Bounding Pair (CBP for short) to be a pair of
two disjoint and non-isotopic homologous curves $\beta$, $\beta'$ on $\Sigma_{g,1}$ such that neither $\beta$ nor $\beta'$ is null-homologous and such that each one bounds a properly embedded disk in $\mathcal{H}_g$. A typical pair is given in Figure \ref{fig intersection}.

A Contractible Bounding  Simple Closed Curve (CBSCC for short) is a non-contractible simple closed curve $\delta$ on $\Sigma_{g,1}$ such that $\Sigma_{g,1} \setminus  \delta$ has two connected components and that bounds a properly embedded disk in $\mathcal{H}_g$. For instance, a curve parallel to the boundary of $\Sigma_{g,1}$ is a CBSCC.

Combining the cited papers of Luft and Johnson we get that if $\beta,\beta'$ is a CBP then the map $T_\beta T_{\beta'}^{-1}$ belongs to $\mathcal{L}_{g,1} \cap \mathcal{TB}_{g,1}$. We call such a map a \emph{CBP-twist}, we also call the intersection group the Luft-Torelli group and we denote it by $\mathcal{LTB}_{g,1}$. In \cite{Joh0}, Johnson proved that opposite twists around Bounding Pairs generate the Torelli group for $g \geq 3$. In this section we prove an analogous theorem for the Luft-Torelli group:

\begin{theorem}\label{thm lufttorelligen}
The Luft-Torelli group $\mathcal{LTB}_{g,1}$ is generated by CBP-twists.
\end{theorem}

\subsection{Reduction to the closed case}\label{subsec reductionclosed}

The reduction to the closed case as many other results in this section are based on the following Lantern Relation, originally  due to Dehn \cite{Deh} and later rediscovered by Johnson \cite{Joh0}. 

\begin{figure}
\input{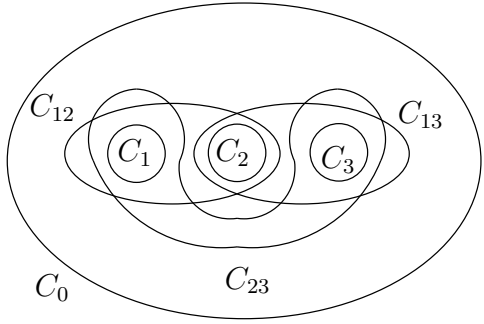}
\caption{Lantern Configuration}\label{fig lanternrel}
\end{figure}

\begin{lemma}[Lantern Relation]\label{lem lanternrel}
Consider a $2$-sphere with $4$ holes (i.e. a lantern). Let the boundary components be $C_0,C_1,C_2,C_3$ and for $1 \leq i < j \leq 3$ denote by $C_{ij}$ a simple curve encircling $C_i$ and $C_j$ (see Figure \ref{fig lanternrel}). Then the following relation  between Dehn twists holds:
\[
T_{C_0} T_{C_1} T_{C_2} T_{C_3} = T_{C_{12}} T_{C_{13}} T_{C_{23}}.
\]
\end{lemma}
Notice that once the four boundary circles are ordered the remaining curves and thus the Lantern Relation are determined.

Recall from Section \ref{sec Heegardsplit} that the kernel of the map $\mathcal{M}_{g,1} \rightarrow \mathcal{M}_{g,\ast}$ is an infinite cyclic group generated by a Dehn twist along a curve $\partial$ parallel to the boundary and that this is a CBSCC. In particular the kernel is contained in $\mathcal{TB}_{g,1}$. Moreover the action of this Dehn twist on the homology of the surface and also on the first homotopy group of the handlebody $\mathcal{H}_g$ is trivial. As a consequence we have a short exact sequence :
\[
\xymatrix{
1 \ar[r] & \Z \ar[r] & \mathcal{LTB}_{g,1} \ar[r] & \mathcal{LTB}_{g,\ast} \ar[r] & 1.
}
\]

The three curves depicted in Figure \ref{fig relbord} plus the boundary curve $\delta$ define a ``Lan\-tern" i.e. a $4$-holed sphere on the surface $\Sigma_g$. Applying the lantern relation of Johnson (see Johnson \cite{Joh0}) one gets~:

\begin{figure}
\input{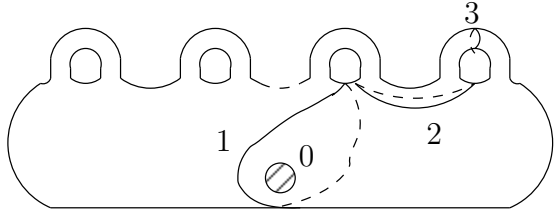}
\caption{Lantern relation for $\delta$}\label{fig relbord}
\end{figure}

\begin{lemma}\label{lem relbord}
The Dehn twist around $\delta$ can be written as a product of CBP-twists.
\end{lemma}

As $\mathcal{LTB}_{g,1}$ is generated by lifts of generators of $\mathcal{LTB}_{g,\ast}$ plus the twist $T_\delta$ and as any CBP-twist in $\mathcal{LTB}_{g,\ast}$ naturally lifts to a CBP-twist in $\mathcal{LTB}_{g,1}$ it is enough to prove 

\begin{proposition}\label{prop lufttoreliiclosed}
The group $\mathcal{LTB}_{g,\ast}$ is generated by CBP-twists.
\end{proposition}

\subsection{Strategy of the proof of Proposition \ref{prop lufttoreliiclosed}}\label{subsec strategy}

Recall from section \ref{subsec actions} that we have two short exact sequences
\[
\xymatrix{
1 \ar[r] & \mathcal{L}_{g,\ast} \ar[r] & \mathcal{B}_{g,\ast} \ar[r] & \text{Aut } \pi_1(\mathcal{H}_{g,\ast}) \ar[r] & 1 \\
1 \ar[r] & \mathcal{TB}_{g,\ast} \ar[r] & \mathcal{B}_{g,\ast} \ar[r] & \mathrm{GL}_{g} (\mathbf{Z}) \ltimes \mathcal{S}_g(\mathbf{Z}) \ar[r] & 1. \\
}
\]

The map $\mathcal{B}_{g,\ast} \rightarrow  \mathrm{GL}_{g} (\mathbf{Z})$ can be identified with the map given by the natural action of $\mathcal{B}_{g,1}$ on the first homology of $\mathcal{H}_{g,\ast}.$ Therefore we have a short exact sequence~:
\[
\xymatrix{
1 \ar[r] & \mathcal{LTB}_{g,\ast} \ar[r] & \mathcal{B}_{g,\ast} \ar[r] & \text{Aut } \pi_1(\mathcal{H}_{g,\ast}) \ltimes \mathcal{S}_g(\mathbf{Z}) \ar[r] & 1. \\
}
\]

In \cite{Niel24a} Nielsen gave an explicit finite presentation  with four generators and 17 relations of the automorphism group of a free group on $g$ generators (see also \cite{MR0422434} Chapter 3), denote this presentation by \[\langle \ x_1, \dots, x_4 \ \vert \ r_1, \dots, r_{17} \ \rangle. \]
The group $\mathcal{S}_g$ is free on $\frac{g(g+1)}{2}$ generators, so we can find a presentation of the form 
\[
\langle \ t_1, \dots, t_{\frac{g(g+1)}{2}} \ \vert \ [t_{i},t_{j}]\ \rangle,
\] where $[t_i,t_j] =t_i^{-1}t_j^{-1}t_it_j$ and $1 \leq i,j \leq \frac{g(g+1)}{2}$. Let the action of $\text{Aut } \pi_1(\mathcal{H}_g)$ on $\mathcal{S}_g(\mathbf{Z})$ be given by expressions of the form: $x_i(t_j) = w_{ij}$ where $w_{ij}$ is a word in the alphabet $t_k$.
Then a presentation of the semi-direct product $\text{Aut } \pi_1(\mathcal{H}_{g,1}) \ltimes \mathcal{S}_g(\mathbf{Z})$ is given by
\[
\langle \ x_1, \dots, x_4, t_1, \dots, t_n \ \vert \ r_1,\dots,r_{17}, [t_{i},t_{j}], x_k^{-1}t_lx_k w_{kl}^{-1} \ \rangle,
\]
where $1 \leq k \leq 4$ and $1 \leq i,j,l \leq \frac{g(g+1)}{2}$.

Assume that we have lifts of the generators $\tilde{x_i}, \tilde{t_j}$ and that these lifts moreover generate $\mathcal{B}_{g,\ast}$. Then a set of normal generators of the group $\mathcal{LTB}_{g,\ast}$ is given by the lifts of the relations $\widetilde{r_i}, [\tilde{t}_{i},\tilde{t}_{j}], \tilde{x}_k^{-1} \tilde{t_l} \tilde{x_k} \tilde{w}_{kl}^{-1} $ as words in the ``lifted" alphabet.

In the mapping class group one has the well-known relation for Dehn twists $\phi T_{\gamma} \phi^{-1} = T_{\phi(\gamma)}$. In particular as the image of a CBP by an element $\phi \in \mathcal{B}_{g,\ast}$ is again a CBP, the group generated by the CBP-twists is normal in $\mathcal{B}_{g,\ast}$. Therefore to prove Theorem \ref{thm lufttorelligen} it is enough to prove

\begin{proposition}\label{prop verifrelateurs}
Under the above hypothesis the lifts \[\widetilde{r_i},[\tilde{t}_{i},\tilde{t}_{j}], \tilde{x}_k^{-1} \tilde{t_l} \tilde{x_k} \tilde{w}_{kl}^{-1}  \in \mathcal{LTB}_{g,\ast} \] are products of CBP-twists.
\end{proposition}

\paragraph*{Proof of Proposition \ref{prop verifrelateurs}}.

It is a classical result of Nielsen \cite{Niel27} that a base-preserving mapping class is determined by its action on the fundamental group of the underlying surface. Geometric descriptions of generators for $\mathcal{B}_{g,\ast}$ were given for instance by Suzuki \cite{Suz} or Luft \cite{MR0500977}. From their work we learn that we need one Dehn twist and 4 particular generators. Here we enlarge the list of Dehn twists to hit each one of the $\frac{g(g+1)}{2}$ generators needed for the group $\mathcal{S}_g(\Z)$.

\subsubsection{Twist generators}\label{subsubsec twistsgen}

\begin{figure}
\input{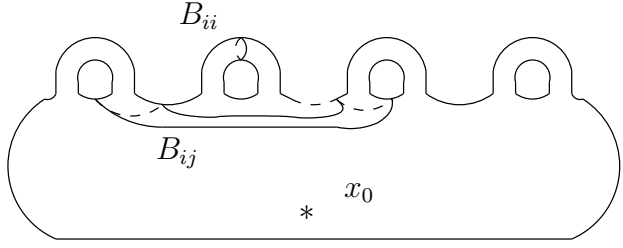}
\caption{Curves for the generating twists}\label{fig twistgenerators}
\end{figure}

We consider the $\frac{g(g+1)}{2}$ curves of Figure \ref{fig twistgenerators}. The curve $B_{ij}$ for $i < j$ goes around the right foot of handles $i$ passes in front of handles $k$ for $i < k <j$ and goes around the left foot of handle $j$, the curve $B_{ii}$ is a meridian of handle $i$.

The corresponding Dehn twists will be denoted respectively $T_{ij}$ and $T_{ii}$. Notice that by construction the homology class of $B_{ij}$ is $b_{ij} = b_i - b_j$ and that of $B_{ii}$ is $b_{ii} = b_i$.

This twists belong to the Twist group $\mathcal{L}_{g,\ast}$ and so act trivially on the homotopy group $\pi_1(\mathcal{H}_g)$. In particular the image  of $T_{ij}$ in $\text{Aut } \pi_1(\mathcal{H}_g) \ltimes \mathcal{S}_{g}(\Z)$ is $(0,t_{b_{ij}})$ where $t_{b_{ij}}$ denotes the transvection along the homology class $b_{ij}$. It is easily verified that these transvections freely generate the group $\mathcal{S}_g(\Z)$.

\subsubsection{Non-twist generators}\label{subsubsec non-twists}

We will keep the names given to these maps by Suzuki in \cite{Suz} but label them according to the generator of $\mathrm{Aut } (\pi_1(\mathcal{H}_g))$ they hit (see \cite[Corollary N1 p.164]{MR0422434}). All elements of the basis of $\pi_1(\Sigma_g)$ that do not appear in the description of the action of a map fixed under the action. We denote by $\sigma_i$ the commutator $\alpha_i^{-1} \beta_i^{-1} \alpha_i \beta_i$.

\begin{enumerate}
\item Cyclic translation of handles, $Q$.
Action on homotopy~:
\[
\begin{array}{rcl}
\alpha_i & \mapsto & \alpha_{i+1} \\
\beta_i & \mapsto & \beta_{i+1}
\end{array}
\]
Indices are counted mod $g$.

\item Twist of knob $1$ $\sigma$.
Action on homotopy~:
\[
\begin{array}{rcl}
\alpha_1 & \mapsto & \alpha_{1}^{-1}\sigma_1^{-1} \\
\beta_1 & \mapsto & \sigma_1 \beta_{1}^{-1}
\end{array}
\]

\item Interchange of knobs $1$ and $2$, $P$.
Action on homotopy~:
\[
\begin{array}{rcl}
\alpha_1 & \mapsto &\sigma_1^{-1} \alpha_2 \sigma_{1} \\
\alpha_2 & \mapsto & \alpha_1 \\
\beta_1 & \mapsto & \sigma_1^{-1} \beta_{2} \sigma_{1} \\
\beta_2 & \mapsto & \beta_1
\end{array}
\]

\item Luft map $U$.
This is a half twist that interchanges the curves $B_{22}$ and $B_{12}$, the boundaries of the two-holed torus which is the support of this map are $B_{11}$ and the curve $C$ depicted in Figure  \ref{fig LuftmapII}.

\begin{figure}
\input{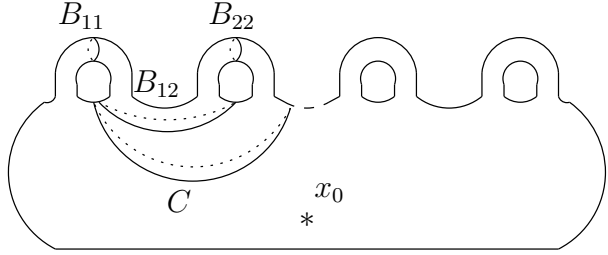}
\caption{Luft map} \label{fig LuftmapII}
\end{figure}

Action on homotopy~:
\[
\begin{array}{rcl}
\alpha_1 & \mapsto & \alpha_1 \alpha_2  \\
\beta_1 & \mapsto & \beta_1 \\
\alpha_2 & \mapsto & \alpha_2^{-1} \beta_2^{-1} \alpha_2^{-1} \beta_2 \alpha_2 \\
\beta_2 & \mapsto & \alpha_2^{-1} \beta_2^{-1} \alpha_1^{-1} \beta_1 \alpha_1 \alpha_2
\end{array}
\]

\end{enumerate}

\subsection{Proof of proposition \ref{prop lufttoreliiclosed}}\label{subsec proof}

According to our strategy of proof we have to lift the relators $r_i, [t_{ij},t_{kl}]$ and $r_i t_{kl} r_i^{-1}w_{ij}^{-1}$ to $\mathcal{B}_{g,1}$ and show that these lifts are products of CBP-twists. We will  deal successively with the twists relators $[t_{ij},t_{kl}]$, the action realtors $x_i t_{kl} x_i^{-1}w_{ij}^{-1}$ and finally with the non-twist relators $r_i$.

Our main tool for recongizing elements that are products of CBP-twists is

\begin{lemma}\label{lem maintool}
Let $\phi \in \mathcal{TB}_{g,\ast}$ be a map. Assume that there exist $g$ disjoint disks $D_i$ properly embedded in the inner handlebody so that $\mathcal{H}_g \setminus  \cup_{i=1}^g D_i$ is a three ball and such that  $\phi(D_i) = D_i$ for $ 1 \leq i \leq g$. Then $\phi$ is a product of CBP-twists. 
\end{lemma}
\begin{proof}
Since $\phi$ acts trivially in homology it can not reverse the orien\-ta\-tions of the boundaries of the discs and therefore we may assume that $\phi$ fixes the disc $D_i$ pointwise. In particular $\phi$ is in the image of the mapping class group relative to the boundary $\mathcal{M}_{0,2g+1}$ of the $2g+1$-holed $3$-ball that is of the complementary of a small neighbourhood of the discs $D_i$. More precisely it is in the kernel of action on homology of this group. This mapping class group is well-known to be isomorphic to the framed pure braid group on $2g$ strands: $\Z^{2g} \times P_{2g}$, where $P_{2g}$ denotes the pure braid group on $2g$ strands.

Without loss of generality we may assume that the discs $D_i$ are our preferred discs $B_{ii}$. Denote the left foot of the $i^{\text{th}}$ handle by $i_0$ and the right foot by $i_1$.

Then the group $\mathcal{M}_{0,2g}$ is generated by the following Dehn twists:
\begin{enumerate}
\item Twists along the curves $A_{i_0 i_0}$ (resp. $A_{i_1 i_1}$) which go to the left (resp.) right foot.
\item Twists along the curves  $A_{i_0 i_1}$  that enclose the two feet of the $i^\text{th}$ handle.
\item Twists the curves  $A_{i_0 j_0}, A_{i_0 j_1}, A_{i_1 j_0}, A_{i_1 j_1}$ for $i< j$ (see Figure \ref{fig twistinPurebraid}). 
\end{enumerate}
\begin{figure}
\input{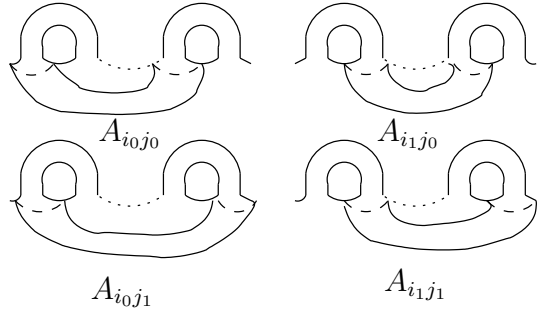}
\caption{Curves for the generators of the Pure Braid group}\label{fig twistinPurebraid} 
\end{figure}

Notice that the twists around $A_{i_0j_1}$ is our preferred Dehn twist $T_{ij}$

If we project compute the action of homology of these twists we get a surjective map $\Z^{2g} \times P_{2g} \rightarrow \mathcal{S}_g(\Z)$ and in terms of the images $\overline{A}_{i_sj_t}$ of the above generators a complete set of generators is given by:
\begin{enumerate}
\item $\overline{A}_{i_0i_1} = 1$, $\overline{A}_{i_0j_1} = \overline{A}_{i_1j_0}$, $\overline{A}_{i_1j_1}= \overline{A}_{i_0j_0}$, $\overline{A}_{i_0j_0} = \overline{A}_{i_0i_0}^2 \overline{A}_{i_1j_0}^{-1} \overline{A}_{j_0j_0}^2$
\item $ [\overline{A}_{i_0j_1}, \overline{A}_{k_0l_1}] = 1$,$ [\overline{A}_{i_0j_1}, \overline{A}_{k_0k_0}] = 1$, $[\overline{A}_{i_0i_0},\overline{A}_{j_0j_0}] = 1$.
\end{enumerate}
Notice that the first series of relations simply reduce the number of relators and the second series is the standard presentation of the free abelian group on the remaining relators.

The kernel of the map $\mathcal{M}_{0,2g} \rightarrow \mathrm{Aut} (H)$ is therefore normally  gene\-rated by the lifts of the above relations that are not relations in $\mathcal{M}_{0,2g}$. We now prove that these lifts as maps in $\mathcal{TB}_{g,1}$ are products of CBP-twists.

There are only three relations in the above list that do not lift obviously to either a relation or a product of CBP-twists.
\begin{enumerate}
\item Relation $\overline{A}_{i_0i_1}$. The twists $A_{i_0i_1}$ are non-trivial mapping classes. For $i=1$, applying the Lantern Relation determined by the four curves in Figure \ref{fig lanternrel2} we can expres $A_{1_0 1_1}$ as a product of CBP-twists. Similar computations hold for the $g-1$ other cases.
\begin{figure}
\input{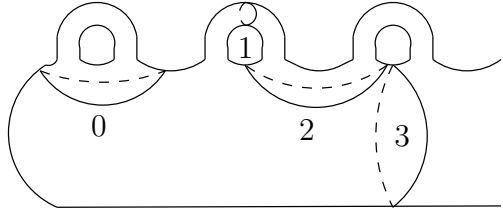}
\caption{Lantern for the Dehn twist $A_{1_0 1_1}$}\label{fig lanternrel2}
\end{figure}
\item  The relation $\overline{A}_{i_0j_0} = \overline{A}_{i_0i_0}^2 \overline{A}_{i_1j_0}^{-1} \overline{A}_{j_0j_0}^2$. We may find a contractible simple closed curve $A'_{i_0j_1}$ that encloses the foots $A_{i_0i_0}, A_{i_1i_0}, A_{j_0j_0}$ defining a Lan\-tern and does not intersect th curves $A_{i_0j_0}$ and $A_{i_1j_0}$. Applying the Lantern relation one finds: 
\[T_{A_{i_0j_0}} T_{A_{i_1j_0}} T_{ A_{i_0 i_0}}^{-2}T_{A_{j_0j_0}}^{-2} = T_{A_{i_0i_1}}^{-1}T_{A'_{i_0j_0}} T_{A_{j_0j_0}}^{-1}T_{A_{i_0i_0}}^{-1}T_{A_{i_1i_1}}.
\]
And this is a product of CBP-twists.
\item The relation $ [\overline{A}_{i_0j_1}, \overline{A}_{k_0l_1}] = 1$ when the curves $A_{i_0j_1}$ and $A_{k_0l_1}$ intersect non-trivially. By definition of the curves this happens if and only if $i< k < j < l.$ 

The relation lifts to
$T_{T_{A_{i_0j_1}(A_{k_0 l_1})}}^{-1} T_{A_{k_0 l_1}}$. As the homology classes of curves $A_{i_0j_1}$ and $A_{k_0 l_1}$ both belong to the Lagrangian $B$, one checks that $T_{A_{i_0 j_1}}(A_{k_0 l_1})$ is homologous to $A_{k_0 l_1}$. In particular the lift is almost a CBP-twists except that the under\-lying curves intersect. It is enough to find a third curve $A'_{k_0l_1}$, disjoint from $A_{i_0 j_1}$ and $A_{k_0 l_1}$, contractible in the inner handlebody and homologous to $A_{k_0 l_1}$ , for then it will be also disjoint from $T_{A_{i_0 j_1}}(A_{k_0 l_1})$ and we will have 
\[
T_{T_{A_{i_0j_1}(A_{k_0 l_1})}}^{-1} T_{A_{k_0 l_1}} = T_{T_{A_{i_0j_1}(A_{k_0 l_1})}}^{-1} T_{A'_{k_0l_1}} T_{A'_{k_0l_1}}^{-1}T_{A_{k_0 l_1}},
\]
a product of CBP-twists. This is done by using curves that go ``through the handles", see Figure \ref{fig throughhangdles} for the case $(i_0,j_1 k_0, l_1) = (1,3,2,4)$.
\begin{figure}
\input{throughhandles.pstex_t}
\caption{Curves for lifting  $[\overline{A}_{1_03_1},\overline{A}_{2_03_1}]$}\label{fig throughhangdles}
\end{figure}
\end{enumerate}
\end{proof}

\subsubsection{Lifts of Twist relators}\label{subsubsec twistsrelators}

One checks that he lifts of the relators $[t_{ij},t_{kl}]$ all leave the curves $B_{ii}$ $1 \leq i \leq g$ invariant and we may apply Lemma \ref{lem maintool}.

\subsubsection{Lifts of action relators}\label{subsubsec actionrelators}

Recall that the generators $t_{ij}$ lift to Dehn twists around the curves $B_{ij}$ and that the lift of the action of the generators $Q, \sigma, P , U$ is conjugation in $\mathcal{B}_{g,\ast}$ by the cooresponding map. For each of the following relations one checks directly that the lifts leave the curves $B_{ii}$ invariant and therefore we may apply in each case Lemma \ref{lem maintool}.

\begin{enumerate}
\item Action of $Q$. The relations to lift are all of the form  $Q (t_{ij}) = t_{i+1j+1}$ (indices mod $g$).

\item Action of $\sigma$. Relations are of the form $\sigma (t_{1i}) = t_{11}^{-2} t_{1j}t_{jj}^{-2} \text{for } 1 < i$, $\sigma (t_{ij}) = t_{ij}$ for $1< i \leq j$ and $\sigma (t_{11}) = t_{11}$.

\item Action of $P$. Relations are $P(t_{11}) = t_{22}$, $P(t_{22}) = t_{11}$, $P(t_{1i}) = t_{2i}$ for $i \geq 3$, $P(t_{2i}) = t_{1i}$ for $i \geq 3$. All other generators are fixed.

\item Action of $U$. Relations are $U(t_{22}) = t_{12}$, $U(t_{12}) = t_{22}$, $U(t_{2i}) = t_{1i} t_{12} t_{11}^{-1}$ for $i \geq 3$. All other generators are fixed.
\end{enumerate}
\subsubsection{Lifts of non-twists relators}\label{subsubsec non-twistsrelators}

Instead of using a case-by-case check we use the following rephrasing of a result of Hirose (see \cite[Theorem $B_\ast$]{MR1643716} and the description of generators therein)~:

\begin{proposition}\label{prop kerhirose}
The kernel of the map $\mathcal{AB}_{g,\ast} \rightarrow \mathrm{Aut }\pi_1(\mathcal{H}_{g})$ is ge\-ne\-ra\-ted by maps which have the following property~:

There exist $g$ properly embedded discs $D_1,\dots D_g$ in $\mathcal{H}_g$  such that
$\mathcal{H}_g \setminus (D_1 \cup \dots D_g)$ is a $3$-ball and such that the map fixes the boundaries of the discs (up to isotopy).
\end{proposition}

In view of Lemma \ref{lem maintool} the maps described in the above Proposition are all products of CBP-twists.

Consider any relation $r$ among the generators of $\mathrm{Aut }\pi_1(\mathcal{H}_{g})$. Since our lifts of the generators of $\mathrm{Aut }\pi_1(\mathcal{H}_{g})$ all belong to $\mathcal{AB}_{g,\ast}$ the lift $\widetilde{r}$ of $r$ belongs to the kernel of the map $\mathcal{AB}_{g,\ast} \rightarrow \mathrm{Aut} \pi_1(\mathcal{H}_g)$. Therefore by the above Proposition \ref{prop kerhirose} and by Lemma \ref{lem maintool}, the lift $\widetilde{r}$ is a product of CBP-twists.

This ends the proof of Proposition \ref{prop verifrelateurs}.

\bibliographystyle{plain}\label{biblography}
\bibliography{3man.bib}

\end{document}